\documentclass[a4paper,10pt]{article}

\usepackage{times}
\usepackage[T1]{fontenc}
\usepackage[latin1]{inputenc}
\usepackage{amsmath,amssymb}
\usepackage{amsthm}
\usepackage{amscd}

\theoremstyle{definition}

\numberwithin{table}{section}
\numberwithin{equation}{section}

\begin{document}
\title{An iteration scheme for monotone operators in Hilbert spaces} 
\author{ Olavi Nevanlinna }
\maketitle

 \begin{center}
{\footnotesize\em 
Aalto University\\
Department of Mathematics and Systems Analysis\\
 email: Olavi.Nevanlinna\symbol{'100}aalto.fi\\[3pt]
}
\end{center}

 \begin{center} { This is a TeX- version of  

{\bf Preprint  No 22, November 1978

Department of Mathematics, University of Oulu, Finland}

with a postscript added at the end}
\end{center}

\begin{abstract}
We give an iteration scheme for finding zeros of maximal monotone operators in Hilbert spaces.  We assume that the operator is defined in the whole space.  The iterates converge strongly to a solution if there exists any,  otherwise they tend to infinity. As an application we  get a strongly convergent minimization scheme for convex functionals in Hilbert spaces.

 \end{abstract}
\bigskip


 
\section{Introduction}

In this note we construct an iteration scheme for finding zeros of monotone operators in Hilbert spaces. The only assumption on the operator needed is that it is defined in the whole space. The scheme produces a sequence which, independently of the initial choice, will converge strongly to a zero of the operator, if there exists any, otherwise it tends to infinity. When applied to  subdifferentials we  obtain a strongly convergent iteration scheme for minimising convex functionals. Recall that the (continuous  version of)  gradient method need not converge strongly, see Baillon [1].

We discussed iteration schemes for monotone operators in an earlier paper [3], but we had to assume that the operators were either bounded or continuous. If the operator is defined in the whole space, then it need not be bounded  but it is always locally bounded [5].  In [4] Petrov and Yurgelas showed how one can in this case find the zero of a strongly monotone operator, Here we point out that their scheme can be modified to cover the general case as well.

 \bigskip
 \section{Iteration scheme}
 
 Let $H$ be a real Hilbert space. A possibly multivalued operator in $H$ is said to be monotone if 
 \begin{equation}
 (y_1-y_2, x_1-x_2) \ge 0 \ \ {\rm for \  all } \ x_i \in D(A) \ {\rm  and }\  y_i  \in Ax_i.
 \end{equation}
 Here we shall assume that $D(A)=H$, and without loss of generality, that $A$ is maximal [2], (otherwise we should call $p\in H$ a (generalised) zero  if $p$ satisfies
 \begin{equation} 
 (y, x-p) \ge 0 \ \ {\rm for \  all } \ x \in H \ {\rm   and }\  y \in Ax,
 \end{equation} 
 which is equivalent to $0 \in \tilde A p$,  where  $\tilde A$ denotes the maximal monotone extension of $A$).

 We now give the algorithm:  Choose any decreasing positive sequence $\{ \theta_\mu\}$ such that the sequence
 $\{\rho_\mu\}$, where  $\rho_\mu = \frac{1}{\theta_\mu+1} -\frac{1}{\theta_\mu}$, is also positive and decreasing and $\lim _{\mu\rightarrow \infty} \theta_\mu =\lim _{\mu\rightarrow \infty} \rho_\mu =0.$

 \smallskip
 {\bf Step 1}   Choose $x_0 \in H$ and set $\mu=0$. Further choose $r_0$ such that for some $y_0 \in Ax_0$
 \begin{equation}
 |x_0 + \frac{1}{\theta_0} y_0| \le r_0.
 \end{equation}

 \smallskip
{\bf Step 2}  If $r_n \le \rho_\mu$, then go to  {Step 3}, otherwise take $y_n \in Ax_n$, set
\begin{equation} 
\lambda_n = \min \{ 1/ 2 \theta_\mu, \ \theta_\mu r^2_n / |y_n + \theta_\mu x_n|^2\}
\end{equation}
and compute
\begin{equation}
x_{n+1} = x_n - \lambda_n (y_n + \theta_\mu x_n)
\end{equation}
and
\begin{equation} 
r_{n+1}^2 = (1-2\theta_\mu \lambda_n)\ r_n^2 + \lambda_n^2 \  | y_n + \theta_\mu x_n|^2.
\end{equation}
Repeat Step 2.

\smallskip
{\bf Step 3 } Redefine $r_n=(1+|z| ) \rho_\mu$, where $z\in A0$.
Increase $\mu$ by $1$ and go back to Step 2. 
 
\bigskip

If we denote the vectors $\theta_\mu x_n$ appearing in Step 2 by $v_n$, then  we have the following
 
 \bigskip

{\bf Theorem}  
  \ {\it Let $A$ be a maximal monotone operator in a Hilbert space $H$ and assume that $D(A)=H$.  Suppose further that $\{x_n\}$ is generated by the iteration scheme described above.  If $A^{-1} \not= \emptyset $, then $x_n$ converges strongly to $p$, where $p$ is the element in $A^{-1} 0$ with minimum norm,  and if $A^{-1}=\emptyset$, then $|x_n| \rightarrow \infty$ in such a way that $v_n$ converges to $ -a^0$, where  $a^0$ is the element in $\overline {R(A)}$ with minimum norm. } 

\bigskip

\section{Proof of Theorem } 

Since $A$ is maximal monotone there exists a unique $p_\mu \in H$ such that
\begin{equation}
 \theta_\mu p_\mu + A p_\mu  \ni 0.
\end{equation}
The algorithm associates to any $n$ an index $\mu= \mu(n) $ by Step 2.   In the following
we denote this index simply by $\mu$.  The convergence of the iterates $x_n$ follows from two facts: first, $|x_n -p_\mu|$ tends to zero as $n$ tends to infinity, and secondly, $p\mu$ tends to minimum solution of  $ Ap \ni 0$, if it exists, or tends to infinity if $0\in Ap$  has no solutions. 

We show first that for all $n\ge 0$
\begin{equation}
|x_n - p_\mu| \le r_n.
\end{equation}
By monotonicity of $A$ we have for all $x \in H$,  $y \in Ax$
\begin{equation}
|x-p_\mu|^2 \le (x-p_\mu+ \frac { 1}{\theta_\mu}(y-q_\mu), \ x-p_\mu) = (x+\frac{1}{\theta_\mu} y, \ x-p_\mu)
\end{equation}
where $q_\mu$ is the vector in $Ap_\mu$ such that  $\theta_\mu p_\mu + q-\mu =0$.  But (3.3)  implies
\begin{equation}
|x-p_\mu | \le |x + \frac{1}{\theta_\mu} y | ,   \ \ {\rm for  \  all   } \ \ x \in H,\  y \in Ax,
\end{equation}
and (3.2) holds for $n=0$ by (2.3).

Suppose now that (3.2) holds for the index $n$, and assume first that $r_n > \rho_\mu$. 
 Subtracting $p_\mu$ from both sides of (2.5),  squaring both sides and using the monotonicity of $A$ yield
 \begin{equation}
| x_{n+1} - p_\mu|^2 \le  (1-2\theta_\mu \lambda_n) |x_n - p_\mu|^2 + \lambda_n^2 |y_n + \theta_\mu x_n|^2.
\end{equation}
 But this together with the definition of $r_{n+1}$ and the assumption $|x_n -p_\mu| \le r_n$
 imply $|x_{n+1} -p_\mu | \le r_{n+1}$.  Suppose then that  $r_n \le \rho_\mu$.  We show that  $|x_n -p_{\mu +1}| \le (1 + |z| ) \rho_\mu$, and hence, after  we redefine $r_n = (1+|z|) \rho_\mu$, we get $|x_n - p_{\mu +1} | \le r_n$ and the previous discussion applies since $r_n >\rho_{\mu+1}$. 
 
 Since $A$ is monotone and $\theta_\mu p_\mu + q_\mu =0$, we obtain
 $$
 |p_\mu -p_{\mu+1}|  \le |p_\mu -p_{\mu+1} + \frac{1}{\theta_{\mu+1}}(q_\mu - q_{\mu+1})|
 = (\frac{\theta_\mu}{\theta_{\mu+1}} - 1) \ | p_\mu |.
 $$
If $z \in A0$, then (3.4) implies $|p_\mu| \le  \frac{1}{\theta_\mu} |z|$  and therefore
$$
|x_n - p_{\mu+1}| \le |x-p_\mu |  + |p_\mu - p_{\mu+1}| \le (1+ |z| ) \rho_\mu
$$
because $ \rho_\mu =  \frac{1}{\theta_{\mu+1}} -  \frac{1}{\theta_{\mu} }$.  Hence (3.2) holds for all $n \ge 0$.

From the definitions of $\lambda_n$ and $r_{n+1}$  se see that, as long as  $r_n > \rho_\mu$ we have  $r_{n+1} < r_n$  and therefore for  all $n$ with  $\mu = \mu(n)  \ge 1$  we have  $r_n \le (1+|z|) \rho_ {\mu-1}$.  But by (3.2), $|x_n-p_\mu | \le (1+|z|) \rho_{\mu-1}$ and  $ |x_n - p_\mu | $ tends to zero as $n\rightarrow \infty$  if $\mu=\mu(n)$ tends to infinity, because  $ \lim_{\mu\rightarrow\infty} \rho_\mu =0$. 

Assume therefore that for some $\mu$ we have $r_n >\rho_\mu $ for all large $n$.  This is a contradiction to a result of Perov and Yurgelas [4], which states that  the iteration (2.4), (2.5), (2.6) (with a fixed $\theta_\mu$ )  converges to  the solution of $\theta_\mu p_\mu + Ap_\mu \ni 0$  and $\lim_{n\rightarrow \infty}=0$.   FOr the convenience of the reader we  outline the proof.  

Consider the iteration (2.4), (2.5), (2.6)  with $\theta_\mu = \theta$ and  with $r_n > \rho_\mu$.  Then $r_n$ is  decreasing and  $\rho=\lim r_n >0 $.  From (2.6) we see that if $\lambda_n = 1/2\theta$, then $r_{n+1} ^2 \le \frac{1}{2} r_n^2$, so that for large enough $n$ we must have  $\lambda_n = \theta r_n^2/ |y_n + \theta x_n|^2$,  and hence for  $n$ large
\begin{equation}
r_{n+1}^2  \le (1-\theta^2 r_n^2/ |y_n + \theta x_n|^2) r_n^2.
\end{equation}
Now $ \lim|y_n + \theta x_n| =\infty$ , since otherwise (3.6) would contradict $r_n \ge \rho>0$.  By [5] $A$  is locally bounded,   and since $|y_n + \theta x_n| \rightarrow \infty$ there exists  an $R>0$ such that $|x_n - p_\mu| >R$ for all large $n$.  Hence $ |x_n-p_\mu|$ is bounded both from above and below.  It follows from Lemma 1 that there exists a constant $\gamma >0$ such that for large $n$
$$
(y_n + \theta x_n, x_n -p_\mu) \ge  \gamma |y_n + \theta x_n|.
$$
This implies  that for large $n$
$$
|x_{n+1}-p_\mu|^2  \le (1-\phi_n) |x_n -p_\mu|^2,
$$
where $\phi_n \sim C/|y_n + \theta x_n|$.  Therefore $\sum 1/ |y_n + x_n| <\infty$ and
$$
\sum |x_{n+1}-x_n| = \sum \lambda_n | y_n + \theta x_n| <\infty.
$$
But if $x_n$ converges, then $ |y_n + \theta x_n$  must be bounded because $A$ is locally bounded  and we arrive into a contradiction.  In the discussion above we used  the following
\bigskip 

{\bf Lemma 1},  [4].   
{\it Let $A$ be a monotone operator in a Hilbert space and assume that $B(x_0,r_0)=\{\xi \ : \ 
|\xi - x_0| \le r_0 \} \subset D(A)$  and $ c_0 =\sup |\eta| < \infty$  for $\eta \in A\xi, \ \xi \in B(x_0, r_0)$.  If $|y| >c_0$ for $y \in Ax$, then
$$
(y, x-x_0) \ge \sigma |y| |x-x_0|,
$$
where $\sigma = \frac{r_0}{r} \{1- (\frac{c_0}{c} )^2\}^{1/2} - \frac{c_0}{c} \{1- (\frac{r_0}{r} )^2\}^{1/2},
$ and $c=|y|, \ r= |x-x_0| $. }

 \bigskip

So far we have shown that  $\lim_{n\rightarrow \infty} |x_n -p_{\mu(n)}| =0$  and $ \mu(n) \rightarrow \infty $ as $n \rightarrow \infty$.  The conclusions of Theorem 2.1 now follow  from Lemma 2, since $p_\mu = J_{1/\theta_\mu}.$ Here $J_\lambda$ denotes the resolvent:  $J_\lambda = (I + \lambda A)^{-1}. $ 

\bigskip

{\bf Lemma 2}, [3].
 {\it Let $A$ be a maximal monotone operator in a Hilbert space and $\lambda>0$. Then
$$
\lim_{\lambda\rightarrow \infty}  \frac{1}{\lambda} J_\lambda = - a^0
$$
where $a^0$ is the element in $\overline {R(A)}$ with minimum norm,  If $A^{-1}0 = \emptyset$  then  $J_\lambda 0  \rightarrow \infty$, and if  $A^{-1}0 \not= \emptyset$ then $J_\lambda 0 \rightarrow p$  as  $\lambda \rightarrow \infty$, where $p $ is the element in $A^{-1}0$ with minimum norm.
}

\bigskip

Let $A_\lambda$ denote the Yosida approximation of $A$:  $A_\lambda = \frac{1}{\lambda} (I-J_\lambda)$.  Then Lemma 2 follows from the identity  $J_\lambda 0 = (A^{-1} )_{1/\lambda} 0$ and from known properties of the resolvent and of the Yosida approximation,  [2].  For a proof see  [3].

 \bigskip
 
 \bigskip

{\bf References}

\bigskip

[1]  J.B. Baillon: Un Exemple Concernant le Comportement Asymptotique de la Solution du Probleme $\frac {du}{dt} + partial \phi (u) \ni 0.$  J. Functional Analysis {\bf 28}, 369-376 (1978)

\smallskip

[2]   H. Brezis:  Operateurs maximaux monotones et semi-groups de contractions dans les espaces de Hilbert.  Math.Studies {\bf 5} , North Holland, 1973
 
 \smallskip
 
 [3]  O. Nevanlinna:  Global Iteration Schemes for Monotone Operators.  Nonlinear Analysis, Theory, Methods \& Applications, Vol. 3, No. 4, pp. 505-514
 
 \smallskip
 
 [4]  A.I. Perov and V.V. Yurgelas:  On the Convergence of an Iterative Process.  Zh.vychil.Mat.mat.Fiz. {\bf 17}, 4, 859-870 (1977). (English translation: U.S.S.R. Comp.Maths Math.Phys.  {\bf 17} , 45-56 (1978)
 
 \smallskip
 
 [5]  R.T.Rockafellar: Local boundedness of nonlinear monotone operators. Michigan Math.J.  {\bf 16}, 397-407
  
  \bigskip
 
  \bigskip

  {\bf  Postscript  in December  2021}
  
  \bigskip

The text above is as it was in the original  report  in November 1978, except that  reference [3]  was then  "to appear".  Here I also included the key  words and MSC classifications.   

I  made  no attempts to improve the, at places clumsy,  language so as  to keep it authentic.  The  library of Oulu University still has a copy of  the original report. 

The note was written while moving  from Madison, Wis.  to Oulu  where I started in September 1978, and as it turned out, the  report was never  published or even submitted.  I soon got  a full professorship in Helsinki University of Technology and was busy running things. 

At times there has been references  to [3],  modifications or generalisations to accretive operators in Banach spaces, in particular if the unit balls are smooth enough,    The work  [4] by Perov and Yurgelas  has  not been referenced much.   As this  note was  partly based on their work,  I find it natural to make the report available and thus also to have a pointer to [4].

 \end{document}